\documentclass[10pt]{article}

\usepackage[english]{babel}
\usepackage{amsmath,amsfonts,amssymb}
\usepackage{graphicx}

\usepackage[a4paper,margin={2cm,2cm}]{geometry}

\newtheorem{thm}{Theorem}[section]
\newtheorem{cor}{Corollary}[section]
\newtheorem{lem}{Lemma}[section]
\newtheorem{prop}{Proposition}[section]

\newenvironment{proof}{\par\noindent{\bf Proof:}}{\hfill$\blacksquare$\par}

\DeclareMathOperator{\sgn}{sgn} % signum

\newcommand{\be}{\begin{equation}}
\newcommand{\ee}{\end{equation}}
\newcommand{\ben}{\begin{equation*}}
\newcommand{\een}{\end{equation*}}
\newcommand{\ba}{\begin{equation}\begin{aligned}}
\newcommand{\ea}{\end{aligned}\end{equation}}

\newcommand{\cD}{\mathcal{D}}

\newcommand{\cF}{\mathcal{F}}
\newcommand{\cG}{\mathcal{G}}

\newcommand{\cT}{\mathcal{T}}
\newcommand{\cS}{\mathcal{S}}

\newcommand{\bI}{\mathbb{I}}

\newcommand{\bR}{\mathbb{R}}

\newcommand{\bN}{\mathbb{N}}

\newcommand{\ec}{\mathrm{e}}
\newcommand{\iu}{\mathrm{i}}

\newcommand{\E}{\mathbf{E}}
\renewcommand{\P}{\mathbf{P}}
\newcommand{\e}{\varepsilon}

\begin{document}

\title{Limit theorems for $p$-variations of
solutions of SDEs driven by additive non-Gaussian stable L\'evy noise}

\author{, Peter Imkeller and Ilya Pavlyukevich}
\author{Claudia\ Hein$^\ast$, Peter\ Imkeller$^\ast$ and Ilya\ Pavlyukevich\footnote{
Postal address:
Institut f\"ur Mathematik,
Humboldt Universit\"at zu Berlin,
Unter den Linden 6,
10099 Berlin, Germany. E--mail:
\texttt{heinc@math.hu-berlin.de} (Claudia Hein),
\texttt{imkeller@math.hu-berlin.de} (Peter Imkeller),
\texttt{pavljuke@math.hu-berlin.de} (Ilya Pavlyukevich).}
}
\date{\today}
% \date{\today\footnote{{\jobname}.tex\hfill \textbf{Preliminary version!!! Do not distribute!!!}}}
\maketitle

\begin{abstract}
In this paper we study the asymptotic properties of the power
variations of stochastic processes of the type $X=Y+L$, where $L$ is an $\alpha$-stable L\'evy process, and $Y$ a perturbation which satisfies some mild Lipschitz continuity assumptions. We establish local functional
limit theorems for the power variation processes of $X$. In case $X$ is a solution of a stochastic differential equation driven by $L$, these limit theorems
provide estimators of the stability index $\alpha$. They are applicable for instance to model fitting problems
for paleo-climatic temperature time series
taken from the Greenland ice core.
\end{abstract}



%\tableofcontents
\numberwithin{equation}{section}
\section{Introduction}

Stochastic differential
equations are being used for quite a while as meso-scopic models for natural phenomena.
In one of their simplest variants, they consist of deterministic differential equations perturbed by a noise term.
The subclass in which the noise is Gaussian arises for instance from microscopic models described by
coupled systems of deterministic differential equations on different time scales, in the limit of infinite fast scale,
as the fluctuations of the slow scale component around its averaged version are considered.
With a view in particular towards the mathematical interpretation of financial time series,
the theory for stochastic differential equations
the noise term of which is given by general (discontinuous and non-Gaussian)
\emph{semimartingales} has received considerable attention during the recent years.\par\smallskip

It is reasonable on a quite general level to model real data by stochastic differential equations.
Usually neither their deterministic nor their noise component can be deduced from first principles
for instance from microscopic models, but may
be selected by statistical inference or \emph{model fit} from the time series they are supposed to interpret.
The central question of the \emph{model selection problem} that motivated this paper asks for
the best choice of the noise term.\par\smallskip

More formally, suppose we wish to model a real time series by the dynamics $X=(X_t)_{t\geq 0}$
of a real valued process of the type
\be
\label{eq:1}
X_t=x+\int_0^t f(s, X_s)\, ds+\eta_t,\quad t\ge 0,
\ee
where the process $\eta$ is the noise component.
Then the problem of a model fit consists in the choice of a
drift term $f$ and a noise term $\eta$, so that the solution of \eqref{eq:1}
is in the best possible agreement with the data of the given time series.\par\smallskip

The example of \eqref{eq:1} which inspired us most
was investigated in the papers \cite{Ditlevsen-99b,Ditlevsen-99a}, where P.\ Ditlevsen
used \eqref{eq:1} as a model fit for temperature data (yearly averages) obtained from the Greenland ice core
describing aspects of the evolution of the Earth's climate
during the last Ice Age, extending over about 100,000 years, in particular the catastrophic warmings and
coolings in the Northern hemisphere, the so-called Dansgaard--Oeschger events
\cite{DansgaardJC-93}. The drift term was chosen as a gradient of a double well
potential (climatic quasi-potential),
the local minima corresponding to the cold and warm meta-stable climate states.
In order to find a good fit for the noise component, P.\ Ditlevsen performed a histogram
analysis for the residuals of the ice core time series, the temperature increments measured between
adjacent data points, i.e.\ years. He conjectured that
the noise may contain a strong $\alpha$--stable component with
$\alpha\approx 1.75$, and plotted an estimate for the
drift term assuming the stationarity of the solution.\par\smallskip

In this paper we resume Ditlevsen's model selection problem for the fit of the noise component
from the perspective of a new testing method to be developed. Following Ditlevsen, we work under
the model assumption that the noise
$\eta$ has an $\alpha$-stable L\'evy component (symmetric or skewed). We search for a test statistics
discriminating well between different $\alpha$, and capable of testing for the right one.
We shall show that this job is well done by the
\emph{equidistant
$p$-variation} --- in the sequel called $p$-variation --- of the process $X$ defined for all $p>0$ as
\ba
V^n_p(X)_t:=\sum_{i=1}^{[nt]}|\Delta_i^n X|^p,
\ea
where $\Delta_i^n X:=X_{\frac{i}{n}}-X_{\frac{i-1}{n}}$ for $1\le i\le n$, $n\geq 1$.
We first observe that for those values of $p$
relevant for our analysis the main contribution to the $p$-variation
of $X$ comes from the $\alpha$-stable component of the noise.
We next prove local limit theorems which hold under very mild assumptions on the
drift term $f$ and allow to determine the stability index
$\alpha$ asymptotically.
We finally use these limit theorems in Section \ref{s:appl} below to analyse the real data
from the Greenland ice core with our methods, and come to an estimate
$\alpha\approx 0.7$, surprisingly quite different from the one obtained by Ditlevsen \cite{Ditlevsen-99b,Ditlevsen-99a}.\par\smallskip

The paper is organised as follows. In section \ref{s:object} we set the stage for stating our main results, which
are applied to the ice core data in section \ref{s:appl}. The remaining sections are devoted to the proof of our
main functional limit theorems, starting with the convergence of the finite dimensional laws in section \ref{fidi}, continuing
with the proof of the tightness of the laws in the Skorokhod topology in section \ref{s:tightness}, and ending with the
robustness proof of the convergence with respect to adding terms of finite variation in section \ref{s:sums}.

In this paper, `$\stackrel{\cD}{\to}$' denotes convergence in the Skorokhod topology, `$\stackrel{d}{\to}$' denotes convergence of finite-dimensional
(marginal) distributions, and `$\stackrel{\text{u.c.p.}}{\to}$' stands for
uniform convergence on compacts in probability. We denote the indicator
function of a set $A$ by $\bI_A$, and $\bar A$ denotes the complement of the set $A$.

\section{Object of study and main results}\label{s:object}

Let $(\Omega,\cF,(\cF_t)_{t\geq 0},\P)$ be a filtered
probability space. We assume that the filtration satisfies the usual
hypotheses in the sense of \cite{Protter-04}, i.e.\ it is right
continuous, and $\cF_0$ contains all the $\P$-null sets of
$\cF$.

For $\alpha \in (0,2]$ let $L=(L_t)_{t\geq 0}$ be an $\alpha$-stable L\'evy process, i.e.\ a
process with right continuous trajectories possessing left side limits (rcll) and stationary independent increments
whose marginal
laws satisfy
\ba
\label{eq:L}
\ln\E e^{\iu\lambda L_t}
&=\begin{cases}
\displaystyle
- t C^\alpha|\lambda|^\alpha
\Big( 1 - \iu  \beta\sgn(\lambda)\tan\frac{\pi\alpha}{2}\Big)
,
& \alpha\neq 1,\\
\displaystyle
 -t C|\lambda|( 1 - \iu\beta \frac{2}{\pi}\sgn(\lambda) \log|\lambda|)  ,
& \alpha= 1,
\end{cases}\quad t\geq 0,\\
\ea
$C>0$ being the scale parameter and $\beta\in[-1,1]$ the skewness.
We adopt the standard notation from
\cite{SamorodnitskyT-94} and write $L_1 \sim S_{\alpha}(C,\beta,0)$.

We also make use of the L\'evy--Khinchin formula for $L$ which takes the
following form in our case (see \cite[Chapter XVII.3]{Feller-71} for details):
\ba
\ln \E e^{\iu\lambda L_1}= \begin{cases}\displaystyle
C_F
\int_{\bR\backslash \{0\}} (e^{\iu\lambda x}-1-\iu\lambda x\bI_{\{|x|<1\}}
\Big[ \frac{1-\beta}{2}  \bI_{\{x<0\}}+  \frac{1+\beta}{2}\bI_{\{x>0\}}\Big]
\frac{dx}{|x|^{1+\alpha}},
& \alpha\neq 1,\\
\displaystyle
C_F
\int_{\bR\backslash \{0\}} (e^{\iu\lambda x}-1-\iu\lambda \sin x)
\Big[ \frac{1-\beta}{2}  \bI_{\{x<0\}}+  \frac{1+\beta}{2}\bI_{\{x>0\}}\Big]
\frac{dx}{|x|^{1+\alpha}}, & \alpha =1,
\end{cases}
\ea
where $C_F$ denotes the scale parameter in Feller's notation and equals
\ba
C_F=\begin{cases}\displaystyle
C^\alpha\Big(\cos\Big(\frac{\pi\alpha}{2}\Big)\Gamma(\alpha)\Big)^{-1},&\alpha\neq 1,\\
\displaystyle C\frac{2}{\pi} ,&\alpha=1.
    \end{cases}
\ea
Recall that a totally asymmetric process $L$ with $\beta=1$ ($\beta=-1$) is called spectrally positive (negative).
A spectrally positive $\alpha$-stable process with $\alpha\in (0,1)$ has a.s.\ increasing sample paths and is called
\textit{subordinator}.

The main results of this paper are presented in the following three theorems.
The first theorem deals with the asymptotic behaviour of the $p$-variation for a stable L\'evy process itself.
As we will see later, this behaviour does not change under perturbations by stochastic processes that satisfy some mild conditions.

\begin{thm}
\label{thm:conv_pvar_levy}
Let $(L_t)_{t\ge 0}$ be an $\alpha$-stable L\'evy process with $L_1 \sim S_{\alpha}(C,\beta,0)$. If $p>\alpha/2$ then
\ba
\big( V_p^n(L)_t - nt\,B_n(\alpha,p)\big)_{t\ge 0} \stackrel{\cD}{\to} (L'_t)_{t\ge 0}\quad \text{as}\ n\to\infty,
\ea
where $L_1'\sim S_{\alpha/p}(C',1,0)$
with the scale
\ba
\label{eq:C}
C'=\begin{cases}\displaystyle
C^p \Big(\frac{\cos( \frac{\pi\alpha}{2p})\Gamma(1-\frac{\alpha}{p})}
{\cos(\frac{\pi\alpha}{2})\Gamma(1-\alpha)} \Big)^{p/\alpha}, & \alpha\neq p,\\
C, & \alpha= p.
   \end{cases}
\ea
 The normalising sequence $(B_n(\alpha, p))_{n\geq 1}$ is deterministic
 and given by
\ba
\label{eq:B}
B_n(\alpha,p)
=\begin{cases}
n^{-p/\alpha} \E |L_1|^p, & p\in(\alpha/2,\alpha),\\
\E \sin\big( n^{-1}|L_1|^\alpha\big), & p=\alpha,\\
0, & p>\alpha.
\end{cases}
\ea
\end{thm}
We remark that the skewness parameter $\beta$ of $L$ does not influence the convergence of $V_p^n(L)_t$ and does not
appear in the limiting process
since the $p$-variation depends only on the absolute values of the increments of $L$. Moreover,
for $p>\alpha$ the limiting process $L'$ is a subordinator.

We next perturb $L$ by some other process $Y$. We impose no restrictions on dependence properties of $Y$ and $L$.
The only conditions on $Y$ concern the behaviour of its $p$-variation.
We formulate two theorems, the first for $p\in(\alpha/2,1)\cup (\alpha,\infty)$, and the second for $p\in (1,\alpha]$.

% , so we have information if $p>\ifrac\alpha 2$.

\begin{thm}
\label{thm:add_pgea}
Let $(L_t)_{t\ge 0}$ be an $\alpha$-stable stochastic process, with $L_1\sim S_{\alpha}(C,\beta,0)$ and $(Y_t)_{t\ge 0}$ be another stochastic process that satisfies
\ba
V_p^n(Y) \stackrel{\text{u.c.p.}}{\to} 0,\quad n\to\infty,
\ea
for some $p\in(\alpha/2,1)\cup (\alpha,\infty)$. Then
\ba
(V_p^n(L+Y)_t - nt\, B_n(\alpha,p))_{t\ge 0} \stackrel{\cD}{\to} (L'_t)_{t\ge 0}\quad\text{ as } n\to\infty,
\ea
with $L'$ and $(B_n(\alpha, p))_{n\geq 1}$ defined in \eqref{eq:C} and \eqref{eq:B}.
\end{thm}

The methods used to prove Theorem~\ref{thm:add_pgea} do not work for $p\in (1,\alpha]$, and in this case we have to impose
stronger conditions on the process $Y$.
% We need Lipschitz-continuity in sets with large probability. Note that the Lipschitz-constant may depend on $\eps$, so this is weaker than Lipschitz-continuity of the process as $K$ may tend to infinity as $\eps$ goes to 0. Nevertheless there are many processes that satisfy this condition, e.g. Lebesgue-integrals with integrands that are bounded in compact subsets of $\set R$.

\begin{thm}\label{thm:add_lip}
Let $(L_t)_{t\ge 0}$ be an $\alpha$-stable stochastic process with $L_1\sim S_{\alpha}(C,\beta,0)$, $\alpha \in (1,2)$
and let $(Y_t)_{t\ge 0}$ be another stochastic process. Let $p\in(1,\alpha]$ and $T>0$.
If $Y$ is such that for every $\delta>0$ there exists $K(\delta)>0$ that satisfies
\ba
\P(|Y_s(\omega)-Y_t(\omega)|\le K(\delta)|s-t| \,\,
\text{\rm for all }\,s,t\in[0,T]) \ge 1-\delta,
\ea
the process $Y$ does not contribute to the limit of $V_p^n(L+Y)$, i.e.\
\ba
(V_p^n(L+Y)_t - nt\, B_n(\alpha,p))_{0\le t\le T}
\stackrel{\cD}{\to} (L'_t)_{0\le t\le T},\quad n\to\infty
\ea
with $L'$ and $(B_n(\alpha, p))_{n\geq 1}$ defined in \eqref{eq:C} and \eqref{eq:B}.
\end{thm}

To be able to study models of the type \eqref{eq:1} we formulate the following
corollary of Theorems~\ref{thm:add_lip} and \ref{thm:add_pgea}
which takes into account that
Lebesgue integral processes are absolutely continuous w.r.t.\ the time variable and thus
qualify as small process perturbations in the sense of the Theorems.

\begin{cor}
\label{cor:sde}
Let $(L_t)_{t\ge 0}$ be an $\alpha$-stable stochastic process, with $L_1\sim S_{\alpha}(C,\beta,0)$, and
let $f:\bR_+\times \bR\to \bR$ be a locally bounded function such that
for some $x\in\bR$ and $T>0$ the unique strong solution
for
\be
\label{eq:dif_X}
X_t=x + \int_0^t f(s,X_s)\, ds+ L_t
\ee
exists on the time interval $[0,T]$.
Then for $p>\alpha/2$ we have
\be
(V_p^n(X)_t - nt\, B_n(\alpha,p))_{0\le t\le T}
\stackrel{\cD}{\to} (L'_t)_{0\le t\le T}
\ee
as $n\to\infty$,
with $L'$ and $(B_n(\alpha, p))_{n\geq 1}$ defined in \eqref{eq:C} and \eqref{eq:B}.
\end{cor}

The functional convergence of power variations of symmetric stable L\'evy processes to stable processes
was first studied by Greenwood in \cite{Greenwood-69}, where
more general \textit{non-equidistant} power variations were considered,
and in particular for $p>\alpha$ the convergence to subordinators was proved.
Further, more general results on power variations of L\'evy processes
are obtained by Greenwood and Fristedt in \cite{GreenwoodF-72}.
In \cite{Jacod-07} and \cite{Jacod-08}, Jacod proves convergence results for $p$-variations
of general L\'evy processes and semimartingales. In particular, several laws of large
numbers and central limit theorems are established. Our results are different from Jacod's
because we consider processes possessing no second moments so that only the
generalised central limit theorem can apply. Moreover, we consider in addition
convergence of perturbed processes.
Corcuera, Nualart and Woerner in \cite{CorcueraNW-07} consider
$p$-variations of a (perturbed) integrated $\alpha$-stable process
of the type $X= Y+\int_0^\cdot u_{s-}\,dL_s$ with some cadlag adapted process $u$.
For $u=1$, our setting results. The paper \cite{CorcueraNW-07} contains a law of large numbers for $0<p<\alpha$ and
a functional central limit theorem for $0<p\leq \alpha/2$.
However, very restrictive conditions on possible perturbation processes $Y$ are
imposed, so that the results are not applicable to processes of the type
\eqref{eq:1}.

\section{Applications to real data\label{s:appl}}

In this section we illustrate our convergence results and show how
they can be applied to the estimation of the stability index $\alpha$.
We emphasise that the conclusions obtained are somewhat
heuristic. Additional work has to be done to provide more precise
statistical properties of the $p$-variation processes as estimators for the stability index,
and to describe the decision rule of the testing procedure along with its quality
features.

We first work with simulated data. Assume they are realizations of the SDE \eqref{eq:dif_X}
where $L$ is a stable process with unknown stability index $\alpha.$
From the data set we extract $m$ samples with $n$ data points each by taking adjacent
non-overlapping groups of $n$ consecutive points. This way we get
the samples $(X^{(i)})_{1\le i\le m}$ with $X^{(i)}=(X^{(i)}_{\frac{1}{n}},X^{(i)}_{\frac{2}{n}},\dots, X^{(i)}_1)$,
$1\le i\le m$. Along each sample we calculate the $p$-variation
\be
V_p^n(X^{(i)})_1 = \sum_{j=1}^n |\Delta^n_j X^{(i)}|^p,\quad 1\le i\le m,
\ee
where the parameter $p$ takes values in some appropriately chosen interval $[p_1, p_2]$.

Due to Corollary \ref{cor:sde}, the random variables $V_p^n(X^{(i)})_1$
converge to the stable random variable \ $L'_1$, possessing stability index $\alpha/p$.
In order to estimate $\alpha$, we compare the law of $V_p^n(X^{(i)})_1$
with some known stable reference law.
Thus for each $p$ we calculate the empirical distribution function of $(V_p^n(X^{(i)})_1)_{1\le i\le n}$ given by
\ba
G_{p,n}(x):=\frac{1}{m}\sum_{i=1}^{m} \bI_{(-\infty, x]}(V_p^n(X^{(i)})_1),\quad x\in\bR,
\ea
and consider for convenience the reference $\frac{1}{2}$-subordinator
with scale parameter $C'>0$ whose probability distribution function $F_{1/2,C'}$ can be explicitly calculated by
\ba
F_{1/2,C'}(x)=\sqrt{\frac{C'}{2\pi}}\int_0^x
\frac{e^{-C'/2y}}{y^{3/2}}\, dy,\quad x\geq 0.
\ea
This will be a candidate for the limiting law of the random variable $L'_1$.
The scale parameter $C'=C'(C)$ of $L'_1$ is connected with the 
scale parameter $C$ of 
$L$ by the relation \eqref{eq:C}. To determine the unknown value of $C$ we  
numerically minimise in $p\in [p_1,p_2]$, $C\in[C_1, C_2]$ the following distance of the Kolmogorov--Smirnov type:
\ba
\label{eq:D}
D_n(C,p)=
\sup_{x\geq 0}|G_{p,n}(x)-F_{1/2,C'(C)}(x)|,
\ea
where $0<p_1<p_2$ and $0<C_1<C_2$ have to be chosen appropriately.
Assume $D_n(C,p)$ attains its unique minimum at $C=C^*$ and $p=p^*$.
Since $V_p^n(X)_1$ converges to a $\frac{1}{2}$-subordinator
if and only if $p=2\alpha$,
we immediately obtain estimates for the scale and stability index of 
the driving process $L$, namely $C^*$ and $\alpha^*=p^*/2$.

\begin{figure}
\begin{center}
\includegraphics[width=6cm]{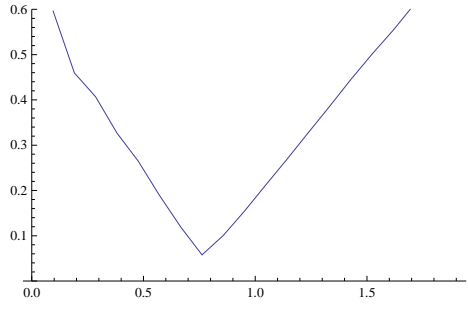}
\end{center}
\caption{$D_n(C^*,p/2)$ for the simulated  data, $L$ is a $0.75$-stable L\'evy process, $n=m=200$\label{fig:kss-s}}
\end{figure}

To test this method, we simulate $m=200$ samples of the data from equation \eqref{eq:dif_X} with $f(\cdot, x)=\cos x$, $x\in\bR$, 
and $L_1\sim S_{0.75}(6.35,0,0)$, $n=200$. We find that the Kolmogorov--Smirnov distance $D_n(C,p)$ attains a unique global
minimum at 
$C^*\approx 6.35$ and 
$p^*\approx 1.5$ corresponding to the true values of $\alpha$ and $C$
(see Fig.~\ref{fig:kss-s}).

We next study the real ice-core data, analysed earlier by Ditlevsen in \cite{Ditlevsen-99b,Ditlevsen-99a}.
The $\log$-calcium signal covers the time period from approximately $90\,150$ to $10\,150$ years before present. We divide it into $m=282$ 
samples, each containing $n=282$ data points. 
Then the Kolmogorov--Smirnov distance is minimised numerically over $p$ and $C$ according to \eqref{eq:D}.

\begin{figure}
\begin{center}
\includegraphics[width=6cm]{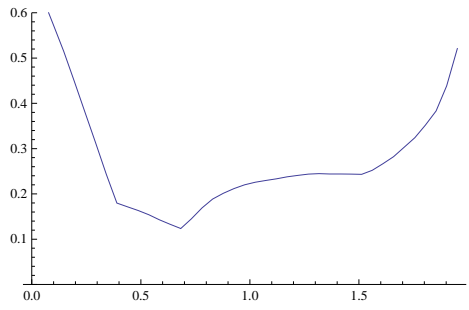}
\end{center}
\caption{$D_n(C^*,p/2)$ for the ice-core data set, $n=m=282$. 
\label{fig:kss-ia}}
\end{figure}

\begin{figure}
\begin{center}
\includegraphics[width=6cm]{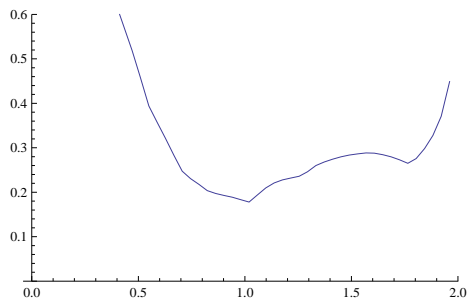}
%\hspace{2cm} \includegraphics[width=6cm]{3dist-icecore-part}
\end{center}
\caption{$D_n(C,p/2)$ for the ice-core data set, $C=3.28$, $n=m=282$. 
% (l.) and for its part from $55\,150$ to $15\,150$ years b.p., $n=m=200$ (r.).
\label{fig:fixedc}}
\end{figure}
%
%
% \begin{figure}
% \begin{center}
% \includegraphics[width=6cm]{dist-icecore-part}
% \end{center}
% \caption{$D(c',p)$ for the ice-core  data (55.000 till 5.000 years b.p.), $n=m=200$\label{fig:kss-ip}}
% \end{figure}

It turns out that $D_n(C,p)$ for the real data also exhibits a 
unique global minimum in the $(C,p)$-domain, which yields the estimate
$\alpha^*\approx 0.7$ for $C^*\approx 7.2$, $D_n(C^*, p^*)\approx 0.1$
(see Fig.~\ref{fig:kss-ia}).
It is striking that our estimate differs from
Ditlevsen's by a quantity very close to 1. This discrepancy can be explained as
follows. It turns out that the function $p\mapsto D_n(C, p)$ has two local minima
for some values of $C$ different from the optimal value $C^*$. For example,
for $C=3.28$ there are two local minima at $\alpha_1\approx 1.02$
and $\alpha_2\approx 1.76$ 
with corresponding
distances $D_n(C,2\alpha_1)\approx 0.178$ and $D_n(C,2\alpha_2)\approx 0.265$. 
Unfortunately the paper \cite{Ditlevsen-99a} only contains the 
estimated value of the stability index $\alpha\approx 1.75$ of the (symmetric)
forcing $L$, and not its scale. It is possible that under some \textit{a priori} 
assumptions on $C$, Ditlevsen's method provides a locally best fit which
is not globally optimal (see Fig.~\ref{fig:fixedc}).

\section{Convergence of the finite dimensional laws of $V^n_p(L)_t$}\label{fidi}
%. Proof of  - theorem \ref{thm:conv_pvar_levy}, part 1}

To prove the convergence of the marginal distributions we use the following theorem
which is a direct result of the well known generalised
central limit theorem
for i.i.d.\ random variables with infinite variance (e.g.\ see Theorem 3 in Feller \cite[Chapter XVII.5]{Feller-71}).

\begin{prop}
\label{thm:clt}
Let $(\eta_i)_{i\geq 1}$ be a sequence of non-negative i.i.d.\
random variables with a regularly varying tail such that
\ba
\P(\eta_1>x)\approx C \frac{2-\alpha}{\alpha} x^{-\alpha/p}\quad\text{as } x\to +\infty.
\ea
for some $\alpha\in (0,2)$, $p>\alpha/2$ and $C>0$.
Then for any $t>0$ we have
\ba
\Big(\frac{t}{n}\Big)^{p/\alpha} \sum_{i=1}^n  \eta_i -
b_{t,n}(\alpha,p)
\stackrel{d}{\to} t^{p/\alpha} Z, \quad\mbox{as}\,\, n\to\infty,
\ea
with
\ba
b_{t,n}(\alpha,p)=
\begin{cases}
\displaystyle
n \Big(\frac{t}{n}\Big)^{p/\alpha}  \E \eta_1,\quad p\in(\alpha/2,\alpha)\\
\displaystyle
n\E\sin \Big(\frac{t\eta_1}{n}\Big),\quad p=\alpha,\\
0,\quad p>\alpha.
\end{cases}
\ea
where $Z\sim S_{\alpha/p} (C',1,0)$ with $C'$ as defined
in \eqref{eq:C}.
\end{prop}
Let $L$ be an $\alpha$-stable L\'evy process as defined in \eqref{eq:L} and
let $p> \alpha/ 2$. To study the finite dimensional distributions
of $V_p^n(L)_t$ we note that
due to the independence of increments of
$L$ it suffices to establish the convergence of marginal laws
for a fixed $t>0$.
Further, the stationarity and independence of increments of $L$ and the self-similarity
property $L_t\stackrel{d}{=}t^{1/\alpha} L_1$ implies that
\be
\label{eq:a}
V_p^n(L)_t =\sum_{i=1}^{[nt]}
|\Delta^n_i L|^p
\stackrel{d}{=}
\sum_{i=1}^{[nt]}
\frac{|\Delta^1_i L|^p}{n^{p/\alpha}},
\ee
with $\Delta^1_i L\stackrel{d}{=}L_1\sim S_\alpha(C,\beta,0)$ being i.i.d.\ random variables.
It is easy to see that the random variables $|\Delta^1_i L|^p$ have a distribution function with a regularly varying tail, namely
\ba
\P ( |\Delta^1_i L|^p>x ) = \P( |L_1| > x^{1/p} )
\approx  C_F\frac{2-\alpha}{\alpha} x^{-\alpha/p}
\quad\text{as }x\to+\infty,
\ea
and thus we can apply Proposition \ref{thm:clt} to the sum \eqref{eq:a}.
Taking into account that $nt/[nt]\to 1$ as $n \to\infty$,
we obtain convergence of the finite dimensional laws in Theorem \ref{thm:conv_pvar_levy}.

\section{Tightness of the laws}\label{s:tightness}

\subsection{Aldous' criterion}

To establish the tightness of the sequence $(V^n_p(L))_{n\geq 1}$ is more complicated.
Although the idea of the proof is based on an application of Aldous' criterion, the technical details
depend strongly on the relationship
between $p$ and $\alpha$.
We first formulate a version of Aldous' criterion for tightness
that is applicable in our case.

On the probability space $(\Omega, \cF,\P)$ we define filtrations
$\cG^n=(\cG^n_t)_{t\geq 0}$ generated by the power variation process
$V_p^n(L)$, i.e.\
\ba
\cG^n_t=\sigma(V_p^n(L)_t:\ s\le t\big)
=\sigma(|L_{\frac{k}{n}}|:\ 0\le k\le [nt]).
\ea
Let $\cT_N^n$ be the set of $\cG^n$-stopping times that are
bounded by $N> 0$.
Then the sequence $V^n_p(L)$ is tight if and only if the following two conditions hold
(see p.\ 350 and p.\ 356 in \cite{JacodS-03}):

\noindent
1.\ for all $\e>0$ and $N>0$ there are $n_0\geq 1$ and $K>0$, so that
for all $n\geq n_0$ we have
\ba
\label{eq:b}
\P(\sup_{t\leq N} V^n_p(L)_t\geq K)\leq \e;
\ea
2.\ for all $N\geq 0$ and $\e>0$ we have
\ba
\label{eq:c}
\lim_{\theta\downarrow 0}\limsup_{n\to \infty}
\sup_{\substack{S,T\in\cT_N^n,\\ S\leq T\leq S+\theta}}
\P(|V^n_p(L)_T-V^n_p(L)_S|\geq \e)=0.
\ea

Condition \eqref{eq:b} is obviously satisfied due to the monotonicity of
the $p$-variation process and
convergence of its marginal distributions to a stable law.

To deal with condition \eqref{eq:c} we note that
since $(V_p^n(L)_t)_{t\ge 0}$ has piecewise constant paths, the filtrations
$\cG^n$ increase only at discrete time instants $t=k/n$, $k\geq 0$, i.e.\
\ba
\cG_t^n = \cG_{\frac{k}{n}}^n\text{ for } \frac{k}{n} \le t< \frac{k+1}n.
\ea
Thus instead of $\cT_N^n$ we can consider in \eqref{eq:c}
a family of
$\cG^n$-stopping times with finitely many values
\ba
\cS_N^n:=\{S\in \cT_N^n : S\,\,\mbox{takes values in}\,\, \{0,\frac{1}{n},\dots,N\}\}.
\ea

\subsection{Tightness for $p>\alpha$}

The case $p>\alpha$ is simple because the compensating sequence $B_{n,t}(\alpha,p)$ vanishes. Thus we can use the monotonicity of the process
$V^n_p(L)_t$.

Let $N\geq 1$ and $\e>0$ be fixed. We use the stationarity and independence of increments of $L$ and convergence the marginal distributions
of $V^n_p(L)$ to $L'$ to obtain the following limit:
\ba
&\lim_{\theta\downarrow 0} \limsup_{n\to\infty}
\sup_{\substack{S,T\in\cS_N^n,\\ S\le T\le S+\theta}} \P \Big( \big| V_p^n(L)_T - V_p^n(L)_S \big| \ge\e \Big)\\
&= \lim_{\theta\downarrow 0} \limsup_{n\to\infty} \sup_{\substack{S,T\in\cS_N^n,\\ S\le T\le S+\theta}} \P \Big( \sum_{i=nS+1}^{nT} |\Delta^n_i L|^p \ge \e \Big)\\
&\leq \lim_{\theta\downarrow 0} \limsup_{n\to\infty}
\sup_{S\in\cS_N^n} \P \Big( \sum_{i=nS+1}^{nS + [n\theta]}
|\Delta^n_i L|^p \ge \e \Big)\\
&= \lim_{\theta\downarrow 0} \limsup_{n\to\infty}
\sup_{S\in\cS_N^n} \sum_{j=0}^{nN} \P
\Big(\sum_{i=j+1}^{j+[n\theta]} |\Delta^n_i L|^p \ge \e, S=\frac jn\Big)\\
&= \lim_{\theta\downarrow 0} \limsup_{n\to\infty} \P
\Big(  \sum_{i=1}^{[n\theta]} |\Delta^n_i L|^p \ge \e\Big)\\
& = \lim_{\theta\downarrow 0} \limsup_{n\to\infty} \P(  V_p^n(L)_\theta \ge \e) \\
&= \lim_{\theta\downarrow 0} \P(  L'_\theta \ge \e)=0.
\ea

\subsection{Tightness for $\frac{\alpha}{2}<p<\alpha$}

The case $p<\alpha$ has to be treated differently since the process
\ba
X_t^n = V_p^n(L)_t -[nt] n^{-p/\alpha} \E|L_1|^p
\ea
need not be monotone.

Let $\theta>0$. As in the preceding section we use the independence and stationarity of increments of $L$ to obtain the estimate
\ba
\label{eq:h}
\sup_{\substack{S,T\in\cS_N^n,\\ S\le T\le S+\theta}}
\P( | X_T^n - X_S^n | \ge\e ) &
\le \sup_{S\in\cS_N^n} \P \Big( \max_{1\le k\le [n\theta]} \Big| \sum_{i=nS+1}^{nS + k} \Big(|\Delta_i^n L|^p - n^{-p/\alpha}
\E|L_1|^p\Big) \Big| \ge \e \Big)\\
& \le \P \Big( \max_{1\le k\le [n\theta]} \Big|\sum_{i=1}^k \Big( |\Delta_i^nL|^p - n^{-p/\alpha} \E|L_1|^p\Big)\Big|  \ge \e\Big)\\
& \le \max_{1\le k\le [n\theta]} 3\P \Big( \Big|\sum_{i=1}^k
\Big( |\Delta_i^nL|^p - n^{-p/\alpha} \E|L_1|^p\Big)\Big|
 \ge \frac\e 3 \Big),
\ea
where the last inequality follows from Lemma 20.2 in \cite{Sato-99}. Now we have to show that the probabilities in the latter inequality
converge to zero uniformly in $k$ as $n\to\infty$ and $\theta\downarrow 0$.

For any $\delta>0$ we can choose $\theta_0>0$ and $n_0\geq 1$ such that
\ba
\label{eq:d}
\P \Big( \Big|\sum_{i=1}^{[n\theta_0]} \Big( |\Delta_i^nL|^p - n^{-p/\alpha} \E|L_1|^p\Big)\Big| \ge \frac \e 6\Big) \le \frac\delta 6
\ea
for all $n\ge \frac{n_0}{2}$. Without loss of generality we can assume that $\theta_0<1$. We will show now that
\ba
\label{eq:e}
\P \Big( \Big|\sum_{i=1}^{[n\theta]} \Big( |\Delta_i^n L|^p - n^{-p/\alpha} \E|L_1|^p\Big)\Big| \ge \frac \e 3\Big) \le \frac\delta 3
\ea
holds for all $\theta<\theta_0$ and $n\ge n_0$. Indeed let us fix $\theta<\theta_0$ and assume first that $\theta<\theta_0/2$.
With the help of the triangle inequality and the estimate \eqref{eq:d}
we can conclude that for $n\ge n_0$
\ba
&\P \Big( \Big|\sum_{i=1}^{[n\theta]} \Big( |\Delta_i^n L|^p - n^{- p/\alpha} \E|L_1|^p\Big)\Big| \ge \frac \e 3\Big)\\
&\le \P \Big( \Big|\sum_{i=1}^{[n\theta_0]} \Big( |\Delta_i^n L|^p
- n^{- p/\alpha} \E|L_1|^p\Big)\Big| \ge \frac \e 6\Big)
+ \P \Big( \Big|\sum_{i=1+[n\theta]}^{[n\theta_0]} \Big( |\Delta_i^n L|^p - n^{-p/\alpha} \E|L_1|^p\Big)\Big| \ge \frac \e 6\Big)\\
&\le \frac\delta 6 + \P \Big( \Big|\sum_{i=1}^{[n\theta_0]-[n\theta]} \Big( |\Delta_i^n L|^p - n^{-p/\alpha} \E|L_1|^p\Big)\Big| \ge \frac \e 6\Big).
\ea
Recalling that $\theta_0-\theta>\theta_0/2$ and $\theta_0<1$ we find $n>n'\ge\frac{n_0}{2}$ such that $[n'\theta_0] = [n\theta_0]-[n\theta]$.
Thus employing the self-similarity of stable processes and \eqref{eq:d} we have
\ba
&\P \Big( \Big|\sum_{i=1}^{[n\theta_0]-[n\theta]} \Big( |\Delta_i^nL|^p - n^{-p/\alpha} \E |L_1|^p \Big)\Big| \ge \frac \e 6\Big)\\
& = \P \Big( \Big|\sum_{i=1}^{[n'\theta_0]} \Big( |\Delta_i^nL|^p - n^{-p/\alpha} \E |L_1|^p \Big)\Big| \ge \frac \e 6\Big)\\
& = \P \Big( \Big(\frac{n'}{n}\Big)^{ p/\alpha}\Big|\sum_{i=1}^{[n'\theta_0]} \Big( |\Delta_i^{n'}L|^p - {(n')}^{-p/\alpha} \E|L_1|^p\Big)\Big| \ge \frac \e 6\Bigg)\\
& \le \P \Big( \Big|\sum_{i=1}^{[n'\theta_0]} \Big( |\Delta_i^{n'}L|^p - (n')^{-p/\alpha} \E|L_1|^p\Big)\Big| \ge \frac \e 6\Big)\leq\frac{\delta}{6}.
\ea
So the inequality \eqref{eq:e} is proved for $\theta<\theta_0/2$ and $n$ big enough. Further, if $\theta_0/2\le\theta<\theta_0$ the proof is even easier because for $n\ge n_0$ we can directly find $n'$ with the same properties as above such that
\ba
\label{eq:f}
&\P \Big( \Big|\sum_{i=1}^{[n\theta]} \Big( |\Delta_i^nL|^p - n^{- p/\alpha} \E |L_1|^p\Big)\Big| \ge \frac \e 3\Big)
= \P \Big( \Big|\sum_{i=1}^{[n'\theta_0]} \Big( |\Delta_i^nL|^p - n^{- p/\alpha} \E|L_1|^p\Big)\Big| \ge \frac \e 3\Big)\\
& \le \P \Big( \Big|\sum_{i=1}^{[n'\theta_0]} \Big( |\Delta_i^{n'}L|^p - (n')^{-p/\alpha} \E|L_1|^p\Big)\Big| \ge \frac \e 3\Big)
 \le \P \Big( \Big|\sum_{i=1}^{[n'\theta_0]} \Big( |\Delta_i^{n'}L|^p
- (n')^{-p/\alpha} \E  |L_1|^p \Big)\Big| \ge \frac \e 6\Big)<
\frac\delta 3.
\ea
Thus we have demonstrated that for any $\delta>0$ and $\e>0$ there
exist $\theta_0>0$ and $n_0\geq 1$ such that
\ba
3\P \Big( \Big|\sum_{i=1}^{[n\theta]} \Big( |\Delta_i^nL|^p - n^{- p/\alpha} \E |L_1|^p \Big)\Big|  \ge \frac\e 3 \Big) \le \delta,\quad n\ge n_0,\ \theta<\theta_0,
\ea
which together with \eqref{eq:h} implies the second condition of Aldous'
criterion, namely
\ba
\lim_{\theta\downarrow 0} \limsup_{n\to\infty} \sup_{\substack{S,T\in\cS_N^n,\\ S\le T\le S+\theta}} \P ( | X_T^n - X_S^n | \ge\e ) \le \delta.
\ea
But $\delta$ was arbitrary.

\subsection{Tightness for $p=\alpha$}

For this part of the proof we need two lemmas that will enable us to control moments of stable processes.

\begin{lem}
\label{thm:ineqmom}
Let $L$ be an $\alpha$-stable L\'evy process and let $k\ge 1$.
Then there are positive numbers $C_{\alpha,k}$, $K_{\alpha,k}$ and $N_\alpha$ such that for $y\ge N_\alpha$ we have
\ba
K_{\alpha,k} y^{(k-1)\alpha} \le \E(|L_1|^{k\alpha}\bI_{\{|L_1|\le y\}})\le C_{\alpha,k} y^{(k-1)\alpha}
\ea
if $k>1$ and
\ba
K_{\alpha,1} \ln y \le \E(|L_1|^{\alpha}\bI_{\{|L_1|\le y\}})
\le C_{\alpha,1} \ln y.
\ea
In particular there is $C_\alpha>0$ such that
$\P(|L_1|\ge y)\le C_\alpha y^{-\alpha}$ for $y\ge N_\alpha$.
\end{lem}

\begin{proof}
We first prove the estimates from above.
Let $F$ denote the probability distribution function of $L_1$. Then there are positive constants $c_\alpha$ and $n_\alpha$ such that for $x\ge n_\alpha$
\ba
\label{eq:p}
1-F(x)\le \frac{c_\alpha}{x^{\alpha}}\quad\text{and}\quad F(-x)\le \frac{c_\alpha}{x^{\alpha}}.
\ea
Integration by parts yields for $y\ge n_\alpha$
\ba
\E\Big(|L_1|^{k\alpha}\bI_{\{|L_1|\le y\}}\Big)
&= \int_{-y}^y |x|^{k\alpha}\, dF(x)\\
% &\qquad = -\int_0^y x^{k\alpha}(F(-x))'\dx - \int_0^y x^{k\alpha} (1-F(x))'\dx
&= - x^{k\alpha}F(-x)\Big|_0^y + \int_0^y k\alpha x^{k\alpha -1}F(-x)\,dx\\
&\qquad\qquad - x^{k\alpha}(1-F(x))\Big|_0^y + \int_0^y k\alpha x^{k\alpha -1}(1-F(x))\,dx\\
&\le 2\int_0^{n_\alpha} k\alpha x^{k\alpha-1}\,dx + 2\int_{n_\alpha}^y k\alpha x^{k\alpha-1}\frac{c_\alpha}{x^\alpha}\,dx\\
& = 2n_\alpha^{k\alpha} + 2c_\alpha\int_{n_\alpha}^y k\alpha x^{(k-1)\alpha-1}\,dx\\
& = 2n_\alpha^{k\alpha} + 2c_\alpha
\begin{cases}
\frac{k}{k-1}\big(y^{(k-1)\alpha}-n_\alpha^{(k-1)\alpha}\big), & k>1,\\
k(\ln y-\ln n_\alpha), & k=1,
\end{cases}\\
& \le 2n_\alpha^{k\alpha} + 2c_\alpha
\begin{cases}
\frac{k}{k-1}y^{(k-1)\alpha}, & k>1,\\ k\ln y, & k=1.
\end{cases}
\ea
Now we can choose constants $C_{\alpha,k}$ and $N_\alpha$ such that the estimates from above are satisfied. The estimates from below follow
analogously from the inequalities
\ba
1-F(x)\ge \frac{k_\alpha}{x^{\alpha}}\quad\text{ and }\quad F(-x)\ge \frac{k_\alpha}{x^{\alpha}}
\ea
for some $k_\alpha>0$ and $x$ big enough.
\end{proof}

\begin{lem}
\label{thm:ineq_cos}
Let $L$ be an $\alpha$-stable L\'evy process. For any $a,b>0$ there exist $C_{\alpha,n}(a,b)$ and $N_\alpha(a,b)$ such that for
$\lambda\in[0,2/a]$ and $n\ge N_\alpha(a,b)$ we have
\ba
\label{eq:g}
\Big| \E \exp\Big({\iu \lambda\frac{|L_1|^\alpha}{n}}\Big) -
\exp\Big(\iu \lambda\E\sin\Big(\frac{|L_1|^\alpha}{n}\Big) \Big)\Big| \le C_{\alpha,n}(a,b)
\ea
with $C_{\alpha,n}(a,b)$ satisfying
\ba
\label{eq:k}
\lim_{b\to 0} \lim_{n\to\infty}C_{\alpha,n}(a,b)nb=0\quad\text{and}\quad \lim_{a\to\infty}\lim_{n\to\infty} C_{\alpha,n}(a,b)n=0.
\ea
\end{lem}

\begin{proof}
We split the left-hand side of \eqref{eq:g} into the real and the imaginary part to obtain the simple estimate
\ba
\Big| \E &\exp\Big({\iu \lambda\frac{|L_1|^\alpha}{n}}\Big) -
\exp\Big(\iu \lambda\E\sin\Big(\frac{|L_1|^\alpha}{n}\Big) \Big)\Big|\\
&\le \Big| \E\cos\Big(\frac{\lambda|L_1|^\alpha}{n}\Big) - \cos\Big(\lambda\E\sin\frac{|L_1|^\alpha}{n}\Big)\Big|
+\Big| \E\sin\Big(\frac{\lambda|L_1|^\alpha}{n}\Big)
- \sin\Big(\lambda\E\sin\frac{|L_1|^\alpha}{n}\Big)\Big|.
\ea
Throughout the proof we shall use the following elementary
inequalities:
\ba
&x-\frac{x^3}{6}\leq \sin x \leq x,\quad x\ge 0,\\
&\cos x \ge 1-\frac{x^2}{2},\quad x\in\bR.
\ea
Let $a,b>0$, $\lambda\in[0,2/a]$ and denote $m_n=m_n(a,b):=na^{1/4} b^{1/2}$.

\noindent
We first estimate the real part.
To this end, we apply Lemma \ref{thm:ineqmom} with $n$ big enough so that
$m_n^{1/\alpha}\ge N_\alpha$ or $n\ge N_\alpha^\alpha a^{-1/4} b^{-1/2}$,
to obtain
\ba
1&\ge\E\cos \Big(\frac{\lambda|L_1|^\alpha}{n}\Big)\\
&\ge \E\cos\Big(\frac{\lambda|L_1|^\alpha}{n}
\bI_{\{|L_1|\le m_n^{1/\alpha}\}}\Big)
- \P(|L_1|\ge m_n^{1/\alpha})\\
&\ge \E\Big(\Big(1- \frac{\lambda^2|L_1|^{2\alpha}}{2n^2}\Big)
\bI_{\{|L_1|\le m_n^{1/\alpha}\}}\Big) - \P(|L_1|\ge m_n^{1/\alpha})\\
&= 1 - \E\Big(\frac{\lambda^2|L_1|^{2\alpha}}{2n^2}
\bI_{\{|L_1|\le m_n^{ 1/\alpha}\}}\Big) - 2\P(|L_1|\ge m_n^{1/\alpha})\\
&\ge 1- \frac{C_{\alpha,2}\lambda^2m_n}{2n^2} - \frac{2C_\alpha}{m_n}\\
&\ge 1- \frac{2C_{\alpha,2}b^{1/2}}{na^{7/4}}
- \frac{2C_\alpha}{na^{1/4} b^{1/2}}=: 1- C_{\alpha,n}^{(1)}(a,b).
\ea
Analogously we get an estimate for the second summand:
\ba
1&\ge \cos\Big( \lambda\E\sin\Big(\frac{|L_1|^\alpha}{n}\Big)\Big)
\ge 1- \frac{\lambda^2}{2}
\Big(\E\sin\Big(\frac{|L_1|^\alpha}{n}\Big)\Big)^2\\
&\ge 1- \frac{\lambda^2}{2} \Big( \E\sin\Big(\frac{|L_1|^\alpha}{n}\bI_{\{|L_1|\le m_n^{1/\alpha}\}} \Big) + \P(|L_1|\ge m_n^{1/\alpha})\Big)^2\\
&\ge 1- \frac{\lambda^2}{2} \Big(
\E\Big(\frac{|L_1|^\alpha}{n}\bI_{\{|L_1|\le m_n^{1/\alpha}\}} \Big)
+ \P(|L_1|\ge m_n^{1/\alpha})\Big)^2\\
&\ge 1- \frac{2}{a^2} \Big( \frac{C_{\alpha,1} \ln m_n}{n\alpha} + \frac{C_\alpha}{m_n}\Big)^2\\
&= 1- \frac{2}{a^2} \Big( \frac{C_{\alpha,1} \ln (na^{1/4} b^{1/2})}{n\alpha} + \frac{C_\alpha}{na^{1/4} b^{1/2}}\Big)^2=: 1- C_{\alpha,n}^{(2)}(a,b).
\ea
Summarising we have
\ba
\label{eq:j}
\Big| \E\cos\Big(\frac{\lambda|L_1|^\alpha}{n}\Big) - \cos\Big(\lambda\E\sin\frac{|L_1|^\alpha}{n}\Big)\Big| \le C_{\alpha,n}^{(1)}(a,b) + C_{\alpha,n}^{(2)}(a,b).
\ea
We now estimate the imaginary part. With analogous arguments we obtain
\ba
&\E\sin\Big(\frac{\lambda|L_1|^\alpha}{n}\Big) - \sin\Big(\lambda\E\sin\Big(\frac{|L_1|^\alpha}{n}\Big)\Big)\\
&\le \E\sin\Big(\frac{\lambda|L_1|^\alpha}{n}\Big)
- \lambda \E\sin\Big(\frac{|L_1|^\alpha}{n}\Big)
+ \frac{\lambda^3}{6} \Big( \E\sin\Big(\frac{|L_1|^\alpha}{n}\Big)\Big)^3\\
&\le \E\Big(\Big(\sin\Big(\frac{\lambda|L_1|^\alpha}{n}\Big) -\lambda\sin\Big(\frac{|L_1|^\alpha}{n}\Big) \Big)
\bI_{\{|L_1|\le m_n^{1/\alpha}\}}\Big)+ (1+\lambda)\P(|L_1|\ge m_n^{1/\alpha}) \\
& +\frac{4}{3a^3}\Big(\frac{C_{\alpha,1} \ln(na^{1/4} b^{1/2})}{n\alpha} + \frac{C_\alpha}{na^{1/4} b^{1/2}}\Big)^3\\
&\le \E\Big(\Big(\frac{\lambda|L_1|^\alpha}{n} -\frac{\lambda|L_1|^\alpha}{n} + \frac{\lambda|L_1|^{3\alpha}}{6n^3}\Big) \bI_{\{|L_1|\le m_n^{1/\alpha}\}}\Big) + \frac{(1+2/a)C_\alpha}{m_n}\\
&+\frac{4}{3a^3}\Big(\frac{C_{\alpha,1} \ln(na^{1/4} b^{1/2})}{n\alpha} + \frac{C_\alpha}{na^{1/4} b^{1/2}}\Big)^3\\
&\le \frac{C_{\alpha,3} b}{3na^{1/2}} + \frac{(1+ 2/a)C_\alpha}{na^{1/4} b^{1/2}} +\frac{4}{3a^3}\Big(\frac{C_{\alpha,1} \ln(na^{1/4} b^{1/2})}{n\alpha} + \frac{C_\alpha}{na^{1/4} b^{1/2}}\Big)^3
 =: C_{\alpha,n}^{(3)}(a,b).
\ea
The lower bound is obtained similarly, and reads
\ba
&\E\sin\Big(\frac{\lambda|L_1|^\alpha}{n}\Big) - \sin\Big(\lambda\E\sin\Big(\frac{|L_1|^\alpha}{n}\Big)\Big)\\
&\ge \E\sin\Big(\frac{\lambda|L_1|^\alpha}{n}\Big)
- \lambda\E\sin\Big(\frac{|L_1|^\alpha}{n}\Big)\\
&\ge \E\Big(\Big(\sin\frac{\lambda|L_1|^\alpha}{n} - \lambda\sin\frac{|L_1|^\alpha}{n}\Big)\bI_{\{|L_1|\le m_n^{1/\alpha}\}}\Big)
 - (1+\lambda)\P( |L_1|\ge m_n^{1/\alpha})\\
&\ge \E\Big(\Big(\frac{\lambda|L_1|^\alpha}{n} - \frac{\lambda^3|L_1|^{3\alpha}}{6n^3} - \frac{\lambda|L_1|^\alpha}{n}\Big)\bI_{\{|L_1|\le m_n^{1/\alpha}\}}\Big) - \frac{(1+\frac{2}{a})C_\alpha}{m_n}\\
&\ge -\frac{4C_{\alpha,3}b}{3a^{5/2}n} - \frac{(1+ 2/a)C_\alpha}{na^{1/4} b^{1/2}}=:C_{\alpha,n}^{(4)}(a,b).
\ea
Combining these estimates with \eqref{eq:j} and denoting
$C_{\alpha,n}(a,b):=C_{\alpha,n}^{(1)}(a,b) +C_{\alpha,n}^{(2)}(a,b) +C_{\alpha,n}^{(3)}(a,b) +C_{\alpha,n}^{(4)}(a,b)$ we obtain inequality
\eqref{eq:g}. Property \eqref{eq:k} of the limits is straightforward.
\end{proof}

Now we can show tightness for $p=\alpha$. We use the inequality analogous to \eqref{eq:h} which states that for any $\e>0$
\ba
&\sup_{\substack{S,T\in\mathcal S_N^n,\\ S\le T\le S+\theta}}
\P \Big( \Big| \sum_{i=nS+1}^{nT} \Big( |\Delta_nL_i|^\alpha - \E \sin\Big(\frac{|L_1|^\alpha}{n} \Big) \Big) \Big|\ge\e \Big)\\
    &\le \max_{1\le k\le [n\theta]}
3\P \Big( \Big| \sum_{i=1}^k \Big( |\Delta_nL_i|^\alpha
- \E \sin\Big(\frac{|L_1|^\alpha}{n} \Big) \Big) \Big|\ge \frac\e 3 \Big).
\ea
Therefore we have to show that
\ba
\lim_{\theta\downarrow 0}\limsup_{n\to\infty} \max_{1\le k\le [n\theta]}
\P ( |X_n^k|\ge\e ) = 0
\ea
with
\ba
X_n^k = \sum_{i=1}^k \Big( |\Delta_nL_i|^\alpha
- \E \sin\Big(\frac{|L_1|^\alpha}{n} \Big) \Big),\quad k,n\geq 1.
\ea
The truncation inequality (see e.g.\ Theorem 5.1 in \cite{Kallenberg-02}) provides an upper bound for the tail of $X_n^k$:
\ba
\label{eq:l}
\P( |X_n^k|\ge\e ) \le \frac\e 2 \int_{-2/\e}^{2/\e} \Big(
1-\E \ec^{\iu \lambda X_n^k}\Big)\,d\lambda
\le 2\sup_{0\le \lambda\le 2/\e} \Big( 1-\Re\E\ec^{\iu \lambda X_n^k}\Big),
\ea
where $\Re c$ denotes the real part of a complex number $c$. Since $X_n^k$ is a sum of i.i.d.\ random variables, its characteristic function can be factorised and we get for $\lambda\ge 0$
\ba
\E\ec^{\iu\lambda X_n^k}
&= \Big[ \E \exp\Big(\iu\lambda\Big(\frac{|L_1|^\alpha}{n} -
\E\sin\Big(\frac{|L_1|^\alpha}{n}\Big)\Big)\Big]^k\\
&= \exp\Big[k\ln \E \exp\Big(\iu\lambda\Big(\frac{|L_1|^\alpha}{n} -
\E\sin\Big(\frac{|L_1|^\alpha}{n}\Big)\Big) \Big]\\
&= \exp\Big[ k\ln
\frac{\E\exp(\iu\lambda \frac{|L_1|^\alpha}{n})}
{\exp(\iu\lambda\E\sin(\frac{|L_1|^\alpha}{n}))}\Big]\\
&= \exp\Big[ k\ln \Big( 1+ \frac{\E\exp(\iu\lambda \frac{|L_1|^\alpha}{n})
- \exp(\iu\lambda\E\sin(\frac{|L_1|^\alpha}{n}))}
{\exp(\iu\lambda\E\sin(\frac{|L_1|^\alpha}{n}))}
\Big)\Big].
\ea
Denote
\ba
\psi_{n,\lambda}:=
 \frac{\E\exp(\iu\lambda \frac{|L_1|^\alpha}{n})
- \exp(\iu\lambda\E\sin(\frac{|L_1|^\alpha}{n}))}
{\exp(\iu\lambda\E\sin(\frac{|L_1|^\alpha}{n}))},
\ea
and note that due to Lemma \ref{thm:ineq_cos} for any $\e,\theta>0$
we can estimate
\ba
|\psi_{n,\lambda}|
= \Big|\E\exp\Big(\iu\lambda \frac{|L_1|^\alpha}{n}\Big)
- \exp\Big(\iu\lambda\E\sin\Big(\frac{|L_1|^\alpha}{n}\Big)\Big) \Big|\le C_{\alpha,n}(\e,\theta),
\ea
with $\lim_{\theta\to 0}\lim_{n\to\infty} C_{\alpha,n}(\e,\theta)\theta n=0$.

Let $\theta>0$ be fixed and let $n_0(\alpha,\e,\theta)$ be such that
$|\psi_{n,\lambda}|\leq 1$ for $n\geq n_0(\alpha,\e,\theta)$.
Recalling the elementary approximations $\ln(1+z)=z+p(z)$ and
$e^z=1+z+q(z)$ with $|p(z)|, |q(z)|\leq |z|^2$, $|z|\leq 1$,
and the estimate \eqref{eq:l} we get
\ba
&\P( |X_n^k|\ge\e)\le
 2\sup_{\lambda\in[0,2/\e]} ( 1-\Re(\exp[ k\ln ( 1+ \psi_{n,\lambda})]) )\\
& = 2\sup_{\lambda\in[0,2/\e]}
\Big( -\Re \Big( k\psi_{n,\lambda} + kp(\psi_{n,\lambda})  + q(k\psi_{n,\lambda} + kp(\psi_{n,\lambda}))\Big) \Big)\\
&\le 2\sup_{\lambda\in[0,2/\e]} \Big( k|\psi_{n,\lambda}| + k|p(\psi_{n,\lambda})| + \big|q(k\psi_{n,\lambda} + kp(\psi_{n,\lambda}))\big| \Big).
\ea
For $n\ge \max(n_0,N_\alpha(\e,\theta))$ with $N_\alpha(\e,\theta)$
as defined in Lemma~\ref{thm:ineq_cos} we already know that
\ba
&k|\psi_{n,\lambda}| + k|p(\psi_{n,\lambda})| \le \theta(n+1) |\psi_{n,\lambda}| + \theta(n+1)|\psi_{n,\lambda}|^2
\le 2\theta(n+1)C_{\alpha,n}(\e,\theta).
\ea
Applying Lemma \ref{thm:ineq_cos} we conclude that for $\theta<\theta_0$, $n\ge \max(n_0,N_\alpha(\e,\theta_0))$ and any $\e>0$
\ba
\lim_{\theta\to 0}\limsup_{n\to\infty} \max_{1\le k\le [n\theta]}
\P( |X_n^k|\ge \e) =0,
\ea
so that condition \eqref{eq:c} of Aldous' criterion holds.

To prove condition \eqref{eq:b} we proceed in exactly the same way as for \eqref{eq:c}, taking advantage of the fact that for fixed $N\geq 1$ the equation
\ba
\lim_{K\to\infty}\lim_{n\to\infty} C_{\alpha,n}(K,N)Nn=0
\ea
holds due to Lemma~\ref{thm:ineq_cos}.

\section{Generalisation to sums of processes}\label{s:sums}

We finally discuss the situation of Theorem \ref{thm:add_lip}, where besides the L\'evy process $L$ another process
$Y$ is given. To see that $V_p^n(L)$ and $V_p^n(L+Y)$ have equivalent asymptotic behaviour we apply Lemma VI.3.31 in \cite{JacodS-03}. Under the conditions of Theorems \ref{thm:add_pgea} and \ref{thm:add_lip} it is enough to show 
that $V_p^n(L+Y)- V_p^n(L)\stackrel{\text{u.c.p.}}{\to} 0$
as $n\to \infty$.

\subsection{Equivalence for $p\le 1$}

Let $L$ be an $\alpha$-stable L\'evy process, $p\in(\alpha/2,1]$ and let $Y$ be such that for $V_p^n(Y)\stackrel{\text{u.c.p.}}{\to} 0$, as
$n\to \infty$. Note that due to the monotonicity properties of $V_p^n(Y)$, the
latter convergence condition is equivalent to
$V_p^n(Y)_N\stackrel{\P}{\to} 0$ for any $N\geq 1$.
Then a simple application of the triangle inequality yields the proof.
In fact, for any $N\geq 1$ we have
\ba
\sup_{0\le t\le N} | V_p^n(L+Y)_t - V_p^n(L)_t| &\le \sum_{i=1}^{nN} \big| |\Delta_i^n(L+Y)|^p - |\Delta_i^nL|^p\big| \\
&\le \sum_{i=1}^{nN} |\Delta_i^n(L+Y) - \Delta_i^nL|^p = V_p^n(Y)_N
\stackrel{\P}{\to} 0,\quad n\to \infty.
\ea

\subsection{Equivalence for $p>\alpha$}
Assume again that
$V_p^n(Y)_t\stackrel{\P}{\to} 0$, $n\to\infty$, for any $t\geq 0$.
Denote $m:=[p]$, so that $p=m+q$ with $q\in [0,1)$. Then for any $N\geq 1$
we have
\ba
&\sup_{t\le N} \left| V_p^n(L+Y)_t - V_p(L)_t \right|
=\sup_{t\le N}  \Big| \sum_{i=1}^{[nt]} |\Delta_i^n(L+Y)|^{m+q} - |\Delta_nL_i|^{m+q} \Big|\\
& = \sup_{t\le N} \Big| \sum_{i=1}^{[nt]} |\Delta_i^n(L+Y)|^q \sum_{k=0}^m
\binom{m}{k}|\Delta_nL_i|^k |\Delta_i^nY|^{m-k} - |\Delta_i^nL|^{m+q} \Big|\\
%&+ \sup_{t\le N} \sum_{i=1}^{[nt]} \Big| |\Delta_i^n(L+Y)|^q |\Delta_i^nL|^m - |\Delta_i^nL|^{m+q}\Big|\\
&\le \sum_{k=0}^{m-1} \binom{m}{k} \sum_{i=1}^{nN} \big| |\Delta_i^n(L+Y)|^q  |\Delta_i^nL|^k |\Delta_i^nY|^{m-k}\Big|\\
& + \sum_{i=1}^{nN} \Big| |\Delta_i^n(L+Y)|^q |\Delta_i^nL|^m - |\Delta_i^nL|^{m+q}\Big|\\
&\le \sum_{k=0}^{m-1} \binom{m}{k} \sum_{i=1}^{nN} |\Delta_i^nL|^{k+q} |\Delta_i^nY|^{m-k} + \sum_{k=0}^{m-1} \binom{m}{k} \sum_{i=1}^{nN} |\Delta_i^nL|^k |\Delta_i^nY|^{m-k+q}\\
&+ \sum_{i=1}^{nN} |\Delta_i^nL|^m\Big| |\Delta_i^n(L+Y)|^q - |\Delta_i^nL|^q\Big|\\
&\le \sum_{k=0}^{m-1} \binom{m}{k} \sum_{i=1}^{nN} |\Delta_i^nL|^{k+q} |\Delta_i^nY|^{m-k}
+ \sum_{k=0}^{m-1} \binom{m}{k} \sum_{i=1}^{nN} |\Delta_i^nL|^k |\Delta_i^nY|^{m-k+q}\\
& + \bI_{\{q>0\}}\sum_{i=1}^{nN} |\Delta_i^nL|^m|\Delta_i^nY|^q .
\ea
The right-hand side of the latter inequality is essentially a sum
of $2m+1$ terms of the type
$\sum_{i=1}^{nN} |\Delta_i^nL|^k |\Delta_i^nY|^{p-k}$, where $k\in[0,p]\cap \bN$.

Applying H\"older's inequality we get
\ba
\sum_{i=1}^{[nt]} |\Delta_nL_i|^k |\Delta_nY_i|^{p-k}
\le \Big( \sum_{i=1}^{[nt]} |\Delta_nL_i|^p \Big)^{k/p}
\Big( \sum_{i=1}^{[nt]} |\Delta_nY_i|^p \Big)^{(p-k)/p}.
\ea
The first factor in the latter formula converges in probability to a finite limit, since $p>\alpha$.
The second factor converges to $0$ in probability due to the assumption
on $Y$, and Theorem \ref{thm:add_pgea} is proven.

\subsection{Equivalence for $\alpha\in (1,2)$, $p\in (1,\alpha]$}

The main technical difficulty of this case arises from the fact that 
the $p$-variation of $L$ for $p<\alpha$ does not exist. In particular,
the events when increments of the stable process $L$ become very large
have to be considered carefully. For $T>0$ and some $c>0$ which will be specified later define the following sets:
\ba
J_c^n(\omega) &:= \{i\in[0,[nT]]:\ |\Delta_i^n L(\omega)|>c\},\\
A_c^n(j) &:= \{\omega\in\Omega:\ |J_c^n(\omega)|=j\}, \quad j=0,\dots, [nT].
\ea
The set $J_c^n$ contains the time instants (in the scale $\frac{1}{n}$), where the
increments of the process $L$ are `large', i.e.\ exceed $c$. The set $A_c^n(j)$ describes
the event that the number of large increments equals $j$.

Let $\delta > 0$, $\e>0$. According to the conditions of Theorem \eqref{thm:add_lip}, let $K':=K(\frac{\delta}{2})>0$, and set
\ba
B=\{\omega: |Y_s(\omega)-Y_t(\omega)|\le K'|s-t| \text{ for all } s,t\in[0,T]\},
\ea
so that $\P(B) \ge 1-\frac{\delta}{2}$.

We estimate
\ba
&\P\Big(\sup_{0\le t\le T}\big|V_p^n(L+Y)_t - V_p^n(L)_t\big|>\e\Big)\\
&
\le \P\Big(\Big\{\sum_{i=1}^{[nT]} \Big||\Delta_i^n(L+Y)|^p - |\Delta_i^nL|^p\Big| > \e\Big\}\cap B\Big)+ \P(\bar{B})\\
& \leq  \sum_{j=0}^{[nT]} \P\Big(\Big\{\sum_{i=1}^{[nT]} \Big||\Delta_i^n(L+Y)|^p - |\Delta_i^nL|^p\Big| > \e \Big\}\cap A_c^n(j)\cap B \Big)+\frac{\delta}{2}\\
&\le \sum_{j=0}^{[nT]} \P\Big(\Big\{\sum_{i\not\in J_c^n}
\Big||\Delta_i^n(L+Y)|^p - |\Delta_i^nL|^p\Big|>\frac{\e}{2}\Big\}\cap A_c^n(j)\cap B\Big)
\\
& + \sum_{j=0}^{[nT]} \P\Big(\Big\{\sum_{i\in J_c^n}
\Big||\Delta_i^n(L+Y)|^p - |\Delta_i^nL|^p\Big|>\frac{\e}{2}\Big\}\cap A_c^n(j)\cap B\Big)+\frac{\delta}{2}\\
&=:D^{(1)}(n,c)+D^{(2)}(n,c)+\frac{\delta}{2}.
\ea
In the following two steps we show that for appropriately chosen $c = c(\e) > 0$
and $n$ big enough, $D^{(1)}(n,c)=0$ and $D^{(2)}(n,c)<\delta/2$. This
will finish the proof.

Step 1. To estimate $D^{(1)}(n,c)$, let $\omega\in B$.
Using the elementary inequality $|a^p-b^p|\le \max\{a,b\}^{p-1}p|a-b|$ which holds for $a,b\ge 0$ and $p\ge 1$ we estimate
\ba
\sum_{i\not\in J_c^n(\omega)} \Big||\Delta_i^n(L(\omega)+Y(\omega))|^p - |\Delta_i^nL(\omega)|^p\Big|
&\leq  \sum_{i\not\in J_c^n(\omega)} p \Big( c+ \frac{K'}{n} \Big)^{p-1}
 \Big||\Delta_i^n(L(\omega)+Y(\omega))| - |\Delta_i^nL(\omega)|\Big|\\
& \le \sum_{i\not\in J_c^n(\omega)} p \Big( c+ \frac{K'}{n} \Big)^{p-1} |\Delta_i^n Y|\\
&\le \sum_{i\not\in J_c^n(\omega)} p \Big( c+ \frac{K'}{n} \Big)^{p-1} \frac{K'}{n}\\
&\le p \bigg( c+ \frac{K'}{n} \bigg)^{p-1} K'T,
\ea
where we have used that the `small' increments of $L$ are bounded by $c$.
So we have
\ba
D^{(1)}(n,c) &= \sum_{j=0}^{[nT]} \P\Big(\Big\{\Big|\sum_{i\not\in J_c^n}|\Delta_i^nL + \Delta_i^nY|^p - |\Delta_i^nL|^p\Big|>\frac\e 2\Big\}\cap A_c^n(j)\cap B\Big)\\
&\le \sum_{j=0}^{[nT]} \P\Big( \Big\{ p\Big(c+\frac{K'}{n}\Big)^{p-1} K'T>\frac{\e}{2}\Big\}\cap A_c^n(j)\cap B \Big)\\
&\le \P\Big( p\Big(c+\frac{K'}{n}\Big)^{p-1} K'T>\frac{\e}{2} \Big)= 0
\ea
for all $n\ge n_1$ with
\ba
n_1 := \Big[ 2K' \Big( \frac{2pK'T}{\e}\Big)^{\frac{1}{p-1}} \Big] \quad\text{and}\quad
c = c(\e) := \frac{1}{4} \Big( \frac{\e}{2pK'T}\Big)^{\frac{1}{p-1}}. 
\ea

Step 2. To estimate $D^{(2)}(n,c)$, let $n\geq \frac{K'}{c}$ so that
for $\omega\in B$ and $i\in J_{c}^n(\omega)$ we have
\ba
\Big| \frac{\Delta_i^nY(\omega)}{\Delta_i^n L(\omega)} \Big|
\le \frac{K'}{nc}\le 1.
\ea
By means of the
elementary inequality $(1+|x|)^p-1\le 3|x|$ which holds for $1\le p\le 2$ and $|x|\le 1$ this implies that
\ba
&\sum_{i\in J_{c}^n(\omega)} \Big||\Delta_i^n(L(\omega)+Y(\omega))|^p
- |\Delta_i^n L(\omega)|^p\Big|\\
&\le \sum_{i\in J_{c}^n(\omega)} |\Delta_i^n L(\omega)|^p \Big[ \Big( 1 + \Big| \frac{\Delta_i^n Y(\omega)}{\Delta_i^n L(\omega)} \Big| \Big)^p-1\Big]\\
&\le \sum_{i\in J_{c}^n(\omega)} |\Delta_i^n L(\omega)|^p\, 3\Big| \frac{\Delta_i^n Y(\omega)}{\Delta_i^n L(\omega)} \Big|\\
&\le \sum_{i\in J_{c}^n(\omega)} |\Delta_i^n L(\omega)|^p \frac{3K'}{nc}.
\ea
This in turn immediately yields
\ba
\label{eq:m}
D^{(2)}(n,c) &\le \sum_{j=1}^{[nT]} \P\Big(A_{c}^n(j) \cap \Big\{\sum_{i\in J_{c}^n}|\Delta_i^n L|^p \frac{3K'}{nc}>\frac{\e}{2}\Big\}\Big)
\le \sum_{j=1}^{[nT]} \P\Big(A_{c}^n(j)\cap  \bigcup_{i\in J_{c}^n} \Big\{|\Delta_i^n L|^p>\frac{\e nc}{6K'j}\Big\} \Big).
\ea
Since all $\Delta_i^n L$, $i=1,\dots,[nT]$, are i.i.d., and only $j$ of them
exceed the threshold $c$, we can estimate the probability for this event precisely.
Indeed, denoting
\ba
p_n:= \P( |\Delta_1^nL| > c ) = \P( |L_1| > c n^{1/\alpha}),
\ea
we continue the estimates in \eqref{eq:m} to get
\ba
\label{eq:n}
&D^{(2)}(n,c) \le \sum_{j=1}^{[nT]} \binom{[nT]}{j}
\P\Big(\Big\{\{|\Delta_i^nL|>c,i=1,\dots, j\}
\cap\{|\Delta_i^nL|\leq c,i=j+1,\dots, [nT]\}\Big\}\cap
\bigcup_{i=1}^j \Big\{ |\Delta_i^nL|^p > \frac{\e nc}{6jK'} \Big\} \Big)
\\
&= \sum_{j=1}^{[nT]} \binom{[nT]}{j} j
\P\Big(   \{|\Delta_i^nL|>c,i=1,\dots, j\}
\cap\{|\Delta_i^nL|\leq c,i=j+1,\dots, [nT]\}\cap
\Big\{ |\Delta_1^nL|^p > \frac{\e nc}{6jK'} \Big\} \Big)
\\
& \le \sum_{j=1}^{[nT]} \binom{[nT]}{j} j p_n^{j-1}(1-p_n)^{[nT]-j}
\P\Big(|\Delta_1^nL|^p > \frac{\e nc}{6jK'}\Big).
\ea
With help of the inequalities \eqref{eq:p}, we obtain the estimates
\ba
\label{eq:pp}
p_n\leq \frac{\tilde c}{c^\alpha n}
\ea
and
\ba
\label{eq:o}
\P\Big(|\Delta_1^nL|^p > \frac{\e nc}{6jK'}\Big)=
\P \Big( |L_1|>\Big(\frac{\e nc}{6K'j}\Big)^{1/p} n^{1/\alpha}\Big)
\le \frac{\tilde c}{n} \Big(\frac{\e nc}{6K'j}\Big)^{-\alpha/p}
\ea
holding for some $\tilde{c}>0$, for all $n\geq n_2 \geq \frac{K'}{c}$ and $1\leq j\leq [nT]$.
Combining \eqref{eq:n} and \eqref{eq:o}, denoting the constant pre-factor by
$C$, and recalling that $\frac{\alpha}{p}\le 2$, we obtain for $n\geq n_2$ that
\ba
D^{(2)}(n,c) \le \frac{C}{p_nn^{1+\alpha/p}}
\sum_{j=1}^{[nT]} \binom{[nT]}{j} j^{1+\alpha/p} p_n^j(1-p_n)^{[nT]-j}
\le \frac{C}{p_n n^{1+\alpha/p}} \sum_{j=1}^{[nT]} \binom{[nT]}{j}
j^3 p_n^j(1-p_n)^{[nT]-j}.
\ea
The sum in the previous formula represents the third moment of a
binomial distribution, and thus can be calculated explicitly.
By means of the asymptotic inequality \eqref{eq:p} we get 
\ba
D^{(2)}(n,c)
&\le C\frac{(nT-2)(nT-1)nT p_n^3 + 3(nT-1)nT p_n^2 + nTp_n}{p_nn^{1+\alpha/p}} \\
&\leq \frac{CT}{n^{\alpha/p}} ( (nTp_n)^2+3nTp_n + 1 ).
\ea
Now choose $n\ge n_3\ge n_2$ big enough to ensure that this expression is smaller than $\frac{\delta}{2}.$ This completes the proof of Theorem \ref{thm:add_lip}.

\par\smallskip

\textsc{Acknowledgements:} P.I.\ and I.P.\ thank DFG SFB 555 \textit{Complex Nonlinear Processes} for financial support. C.H.\ thanks DFG IRTG 
\textit{Stochastic Models of Complex Processes}
for financial support.
C.H.\ is grateful to R.\ Schilling for his valuable comments. The authors
thank P.\ Ditlevsen for providing the ice-core data.

    \bibliographystyle{plain}
%   \bibliography{ref}

    \bibliography{biblio-new}

\end{document}